# Utility of Choice: An Information Theoretic Approach to Investment Decision-making

M. Khoshnevisan
Griffith University
Gold Coast, Queensland Australia

Sukanto Bhattacharya
Bond University
Gold Coast, Queensland Australia

Florentin Smarandache
University of New Mexico - Gallup, USA

**Abstract:**
In this paper we have devised an alternative methodological approach for quantifying utility in terms of expected information content of the decision-maker's choice set. We have proposed an extension to the concept of utility by incorporating extrinsic utility; which we have defined as the utility derived from the element of choice afforded to the decision-maker by the availability of an object within his or her object set. We have subsequently applied this extended utility concept to the case of investor utility derived from a structured, financial product – an custom-made investment portfolio incorporating an endogenous capital-guarantee through inclusion of cash as a risk-free asset, based on the Black-Scholes derivative-pricing formulation. We have also provided instances of potential application of information and coding theory in the realms of financial decision-making with such structured portfolios, in terms of transmission of product information.

**Introduction:**
In early nineteenth century most economists conceptualized utility as a psychic reality – cardinally measurable in terms of *utils* like distance in kilometers or temperature in degrees centigrade. In the later part of nineteenth century Vilfredo Pareto discovered that all the important aspects of demand theory could be analyzed ordinally using geometric devices, which later came to be known as "indifference curves". The indifference curve approach effectively did away with the notion of a cardinally measurable utility and went on to form the methodological cornerstone of modern microeconomic theory.

An indifference curve for a two-commodity model is mathematically defined as the locus of all such points in $E^2$ where different combinations of the two commodities give the same level of satisfaction to the consumer so as the consumer is indifferent to any particular combination. Such indifference curves are always *convex to the origin* because of the operation of the *law of substitution*. This law states that the scarcer a commodity becomes, the greater becomes its relative substitution value so that its marginal utility rises relative to the marginal utility of the other commodity that has become comparatively plentiful.

In terms of the indifference curves approach, the problem of utility maximization for an individual consumer may be expressed as a constrained non-linear programming problem that may be written in its general form for an n-commodity model as follows:

$$\text{Maximize } U = U(C_1, C_2 \ldots C_n)$$
$$\text{Subject to } \Sigma C_j P_j \leq B$$
$$\text{and } C_j \geq 0, \text{ for } j = 1, 2 \ldots n \tag{1}$$

If the above problem is formulated with a strict equality constraint i.e. if the consumer is allowed to use up the entire budget on the n commodities, then the utility maximizing condition of consumer's equilibrium is derived as the following first-order condition:

$$\partial U/\partial C_j = (\partial U/\partial C_j) - \lambda P_j = 0 \text{ i.e.}$$
$$(\partial U/\partial C_j)/P_j = \lambda^* = \text{constant, for } j = 1, 2 \ldots n \tag{2}$$

This pertains to the classical economic theory that in order to maximize utility, individual consumers necessarily must allocate their budget so as to *equalize the ratio of marginal* *utility to price* for every commodity under consideration, with this ratio being found equal to the optimal value of the *Lagrangian multiplier* $\lambda^*$.

However a rather necessary pre-condition for the above indifference curve approach to work is $(U_{C_1}, U_{C_2} \ldots U_{C_n}) > 0$ i.e. the marginal utilities derived by the consumer from each of the n commodities must be positive. Otherwise of course the problem degenerates. To prevent this from happening one needs to strictly adhere to the law of substitution under all circumstances. This however, at times, could become an untenable proposition if measure of utility is strictly restricted to an intrinsic one. This is because, for the required condition to hold, each of the n commodities necessarily *must always have a positive intrinsic utility* for the consumer. However, this would invariably lead to anomalous reasoning like the intrinsic utility of a woolen jacket being independent of the temperature or the intrinsic utility of an umbrella being independent of rainfall.

Choice among alternative courses of action consist of trade-offs that confound subjective probabilities and marginal utilities and are almost always too coarse to allow for a meaningful separation of the two. From the viewpoint of a classical statistical decision theory like that of *Bayesian inference* for example, failure to obtain a correct representation of the underlying behavioral basis would be considered a major pitfall in the aforementioned analytical framework.

Choices among alternative courses of action are largely determined by the relative degrees of belief an individual attaches to the prevailing uncertainties. Following Vroom (Vroom; 1964), the *motivational strength* $S_n$ of choice $c_n$ among N alternative available choices from the choice set C = $\{c_1, c_2 \ldots c_N\}$ may be ranked with respect to the multiplicative product of the relative reward r ($c_n$) that the individual attaches to the consequences resulting from the choice $c_n$, the likelihood that the choice set under consideration will yield a positive intrinsic utility and the respective probabilities p{r ($c_n$)} associated with r ($c_n$) such that:

$$S_{max} = \text{Max}_n [r (c_n) \times p (U_{r(C)} > 0) \times p\{r (c_n)\}], n = 1, 2 \ldots N \qquad (3)$$

Assuming for the time-being that the individual is calibrated with perfect certainty with respect to the intrinsic utility resulting from a choice set such that we have the condition $p (U_{r(C)} > 0) = \{0, 1\}$, the above model can be reduced as follows:

$$S_{max} = \text{Max}_k [r(c_k) \times p\{r(c_k)\}], k = 1, 2 \ldots K \text{ such that } K < N \tag{4}$$

Therefore, choice A, which entails a large reward with a low probability of the reward being actualized could theoretically yield the same motivational strength as choice B, which entails a smaller reward with a higher probability of the reward being actualized.

However, we recognize the fact that the *information conveyed* to the decision-maker by the outcomes would be quite different for A and B though their values may have the same mathematical expectation. Therefore, whereas intrinsic utility could explain the ranking with respect to *expected value* of the outcomes, there really has to be another dimension to utility whereby the **expected information** is considered – that of *extrinsic utility*. So, though there is a very low probability of having an unusually cold day in summer, the information conveyed to the likely buyer of a woolen jacket by occurrence of such an aberration in the weather pattern would be quite substantial, thereby validating a extended substitution law based on an *expected information measure* of utility. The specific objective of this paper is to formulate a mathematically sound theoretical edifice for the formal induction of extrinsic utility into the folds of statistical decision theory.

**A few essential working definitions**

*Object:* Something with respect to which an individual may perform a specific goal-oriented behavior

*Object set:* The set O of a number of different objects available to an individual at any particular point in space and time with respect to achieving a goal where $n\{O\} = K$

**Choice:** A path towards the sought goal emanating from a particular course of action - for a single available object within the individual's object set, there are two available choices - either the individual takes that object or he or she does not take that object. Therefore, generalizing for an object set with K alternative objects, there can be $2^K$ alternative courses of action for the individual

**Choice set:** The set C of all available choices where $C = \boldsymbol{P} O$, $n\{C\} = 2^K$

**Outcome:** The relative reward resulting from making a particular choice

Decision-making is nothing but goal-oriented behavior. According to the celebrated *theory of reasoned action* (Fishbain; 1979), the immediate determinant of human behavior is the intention to perform (or not to perform) the behavior. For example, the simplest way to determine whether an individual will invest in Acme Inc. equity shares is to ask whether he or she *intends* to do so. This does not necessarily mean that there will always be a perfect relationship between intension and behavior. However, there is no denying the fact that people usually tend to act in accordance with their intensions.

However, though intention may be shaped by a positive intrinsic utility expected to be derived from the outcome of a decision, the *ability* of the individual to actually act according to his or her intention also needs to be considered. For example, if an investor truly intends to buy a call option on the equity stock of Acme Inc. even then his or her intention cannot get translated into behavior if there is no exchange-traded call option available on that equity stock. Thus we may view the *additional element of choice* as a measure of extrinsic utility. ***Utility is not only to be measured by the intrinsic want-satisfying capacity of a commodity for an intending individual but also** by the availability of the particular commodity at that point in space and time to enable that individual to act according to his or her intension*. Going back to our woolen jacket example, though the intrinsic utility of such a garment in summer is practically zero, the extrinsic utility afforded by its mere availability can nevertheless suffice to uphold the law of substitution.

**Utility and thermodynamics**

In our present paper we have attempted to extend the classical utility theory applying the entropy measure of information (Shannon, 1948), which by itself bears a direct constructional analogy to the **Boltzmann equation** in thermodynamics. There is some uniformity in views among economists as well as physicists that a functional correspondence exists between the formalisms of economic theory and classical thermodynamics. The laws of thermodynamics can be intuitively interpreted in an economic context and the correspondences do show that thermodynamic entropy and economic utility are related concepts sharing the same formal framework. Utility is said to arise from that component of thermodynamic entropy whose change is due to irreversible transformations. This is the standard *Carnot entropy* given by **dS = $\delta$Q/T** where S is the entropy measure, Q is the thermal energy of state transformation (irreversible) and T is the absolute temperature. In this paper however we will keep to the information theoretic definition of entropy rather than the purely thermodynamic one.

**Underlying premises of our extrinsic utility model**

1. Utility derived from making a choice can be distinctly categorized into two forms:

    (a) Intrinsic utility ($U_{r(C)}$) – the intrinsic, non-quantifiable capacity of the potential outcome from a particular choice set to satisfy a particular human want under given circumstances; in terms of expected utility theory $U_{r(C)} = \Sigma\ r(c_j)\ p\{r(c_j)\}$, where j = 1, 2 … K and

    (b) Extrinsic utility ($U_X$) – the additional possible choices afforded by the mere availability of a specific object within the object set of the individual

2. An choice set with n (C) = 1 (i.e. when K = 0) with respect to a particular individual corresponds to lowest (zero) extrinsic utility; so $U_X$ cannot be negative

3. The law of diminishing marginal utility tends to hold in case of Ux when an individual repeatedly keeps making the same choice to the exclusion of other available choices within his or her choice set

Expressing the frequency of alternative choices in terms of the probability of getting an outcome $r_j$ by making a choice $c_j$, the generalized extrinsic utility function can be framed as a modified version of *Shannon's entropy function* as follows:

$$U_X = -K\ \Sigma_j\ p\{r(c_j)\}\ \log_2 p\{r(c_j)\},\ j = 1, 2 \ldots 2^K \quad (5)$$

The multiplier -K = -n (O) is a scale factor somewhat analogous to the **Boltzmann constant** in classical thermodynamics with a reversed sign. Therefore general extrinsic utility maximization reduces to the following non-linear programming problem:

$$\text{Maximize } U_X = -K\ \Sigma_j\ p\{r(c_j)\}\ \log_2 p\{r(c_j)\}$$
$$\text{Subject to } \Sigma\ p\{r(c_j)\} = 1,$$
$$p\{r(c_j)\} \geq 0;\ \text{and}$$
$$j = 1, 2 \ldots 2^K \quad (6)$$

Putting the objective function into the usual Lagrangian multiplier form, we get

$$Z = -K \Sigma p\{r(c_j)\} \log_2 p\{r(c_j)\} + \lambda (\Sigma p\{r(c_j)\} - 1) \qquad (7)$$

Now, as per the first-order condition for maximization, we have

$$\partial Z/\partial p\{r(c_j)\} = -K(\log_2 p\{r(c_j)\} + 1) + \lambda = 0 \text{ i.e.}$$

$$\log_2 p\{r(c_j)\} = \lambda/K - 1 \qquad (8)$$

Therefore; for a pre-defined K; $p\{r(c_j)\}$ is independent of j, i.e. all the probabilities are necessarily equalized to the constant value **$p\{r(c_j)\}^* = 2^{-K}$** at the point of maximum $U_X$.

It is also intuitively obvious that when $p\{r(c_j)\} = 2^{-K}$ for $j = 1, 2, \ldots 2^K$, the individual has the *maximum element of choice in terms of the available objects within his or her object set*. For a choice set with a single available choice, the extrinsic utility function will be simply given as **$U_X = -p\{r(c)\} \log_2 p\{r(c)\} - (1 - p\{r(c)\}) \log_2 (1 - p\{r(c)\})$**. Then the slope of the marginal extrinsic utility curve will as usual be given by **$d^2U_X/dp\{r(c)\}^2 < 0$**, and this can additionally serve as an alternative basis for intuitively deriving the generalized, downward-sloping demand curve and is thus a valuable theoretical spin-off!

Therefore, though the mathematical expectation of a reward resulting from two mutually exclusive choices may be the same thereby giving them equal rank in terms of the intrinsic utility of the expected reward, the expected information content of the outcome from the two choices will be quite different given different probabilities of getting the relative rewards. The following vector will then give a *composite measure of total expected utility* from the object set:

$$U = [U_r, U_X] = [\Sigma r(c_j) p\{r(c_j)\}, -K \Sigma_j p\{r(c_j)\} \log_2 p\{r(c_j)\}], j = 1, 2 \ldots 2^K \qquad (9)$$

Now, having established the essential premise of formulating an extrinsic utility measure, we can proceed to let go of the assumption that an individual is calibrated with perfect certainty about the intrinsic utility resulting from the given choice set so that we now look at the full Vroom model rather than the reduced version. If we remove the restraining condition that **p (U$_{r\ (C)}$ > 0) = {0, 1}** and instead we have the more general case of **0 ≤ p (U$_{r(C)}$ > 0) ≤ 1**, then we introduce another probabilistic dimension to our choice set whereby the individual is no longer certain about the nature of the impact the outcomes emanating from a specific choice will have on his intrinsic utility. This can be intuitively interpreted in terms of the *likely opportunity cost* of making a choice from within a given choice set to the exclusion of all other possible choice sets. For the particular choice set C, if the likely opportunity cost is less than the potential reward obtainable, then **U$_{r\ (c)}$ > 0**, if opportunity cost is equal to the potential reward obtainable, then **U$_{r(C)}$ = 0**, else if the opportunity cost is greater than the potential reward obtainable then **U$_{r\ (C)}$ < 0**.

Writing **U$_{r(C)}$ = Σ$_j$ r (c$_j$) p{r (c$_j$)}**, j = 1, 2 … N, the *total expected utility vector* now becomes:

**[U$_{r(C)}$, U$_X$] = [Σ$_j$ r (c$_j$) p{r (c$_j$)}, - K Σ p {r (c$_j$)| U$_{r(C)}$ > 0} log$_2$ p {r (c$_j$)| U$_{r(C)}$ > 0}]**, j = 1, 2 … N          **(10)**

Here **p {r (c$_j$)| U$_{r(C)}$ > 0}** may be estimated by the standard ***Bayes criterion*** as under:

**p {r (c$_j$)| U$_{r(c)}$ >0} = [p {(U$_{r(C)}$ ≥0|r (c$_j$)} p {(r (c$_j$)}][Σ$_j$ p {(U$_{r(C)}$ >0|r (c$_j$)} p {(r (c$_j$)}]$^{-1}$**     **(11)**

**A practical application in the realms of Behavioral Finance - Evaluating an investor's extrinsic utility from capital-guaranteed, structured financial products**

Let a structured financial product be made up of a basket of n different assets such that the investor has the right to claim the return on the best-performing asset out of that basket after a stipulated holding period. Then, if one of the n assets in the basket is the risk-free asset then the investor gets assured of a minimum return equal to the risk-free rate i on his invested capital at the termination of the stipulated holding period. This effectively means that his or her investment becomes endogenously capital-guaranteed as the terminal wealth, even at its worst, cannot be lower in value to the initial wealth plus the return earned on the risk-free asset minus a finite cost of **portfolio insurance**.

Therefore, with respect to each risky asset, we can have a *binary response from the investor in terms of his or her funds-allocation decision* whereby the investor either takes funds out of an asset or puts funds into an asset. Since the overall portfolio has to be *self-financing* in order to pertain to a **Black-Scholes** kind of pricing model, funds added to one asset will also mean same amount of funds removed from one or more of the other assets in that basket. If the basket consists of a single risky asset s (and of course cash as the risk-free asset) then, if $\eta_s$ is the amount of re-allocation effected each time with respect to the risky asset s, the two alternative, mutually exclusive choices open to the investor with respect to the risky asset s are as follows:

(1) **C ($\eta_s \geq 0$)** (funds left in asset s), with associated outcome **r ($\eta_s \geq 0$)**; and

(2) **C ($\eta_s < 0$)** (funds removed from asset s), with associated outcome **r ($\eta_s < 0$)**

Therefore what the different assets are giving to the investor apart from their intrinsic utility in the form of higher expected terminal reward is some *extrinsic utility in the form of available re-allocation options*. Then the expected present value of the final return is given as follows:

$$E(r) = \text{Max}\,[w,\ \text{Max}_j\{e^{-it}\,E(r_j)_t\}],\ j = 1, 2 \ldots 2^{n-1} \qquad (12)$$

In the above equation i is the rate of return on the risk-free asset and t is the length of the investment horizon in continuous time and w is the initial wealth invested i.e. ignoring insurance cost, if the risk-free asset outperforms all other assets **$E(r) = we^{it}/e^{it} = w$**.

Now what is the probability of each of the (n – 1) risky assets performing worse than the risk-free asset? Even if we assume that there are some cross-correlations present among the (n – 1) risky assets, given the statistical nature of the risk-return trade-off the joint probability of these assets performing worse than the risk-free asset will be very low over moderately long investment horizons. And this probability will keep going down with every additional risky asset added to the basket. Thus each additional asset will *empower the investor* with additional choices with regards to re-allocating his or her funds among the different assets according to their observed performances.

Intuitively we can make out that the extrinsic utility to the investor is indeed maximized when there is *an equal positive probability of actualizing each outcome $r_j$ resulting from $\eta_j$ given that the intrinsic utility $U_{r(C)}$ is greater than zero.* By a purely economic rationale, *each additional asset introduced into the basket will be so introduced if and only if it significantly raises the expected monetary value of the potential terminal reward.* As already demonstrated, the extrinsic utility maximizing criterion will be given as under:

$$p(r_j | U_{r(C)} > 0)^* = 2^{-(n-1)} \text{ for } j = 1, 2 \ldots 2^{n-1} \qquad (13)$$

The *composite utility vector from the multi-asset structured product* will be as follows:

$$[U_{r(C)}, U_X] = [E(r), -(n-1)\Sigma\, p\{r_j | U_{r(C)} > 0\} \log_2 p\{r_j | U_{r(C)} > 0\}], j = 1, 2 \ldots 2^{n-1} \qquad (14)$$

Choice set with a structured product having two risky assets (and cash):

| 0 | 0 |
|---|---|
| 1 | 0 |
| 0 | 1 |
| 1 | 1 |

That is, the investor can remove all funds from the two risky assets and convert it to cash (the risk-free asset), or the investor can take funds out of asset 2 and put it in asset 1, or the investor can take funds out of asset 1 and put it in asset 2, or the investor can convert some cash into funds and put it in both the risky assets. Thus there are 4 alternative choices for the investor when it comes to re-balancing his portfolio.

Choice set with a structured product having three risky assets (and cash):

| 0 | 0 | 0 |
|---|---|---|
| 0 | 0 | 1 |
| 0 | 1 | 0 |
| 0 | 1 | 1 |
| 1 | 0 | 0 |
| 1 | 0 | 1 |
| 1 | 1 | 0 |
| 1 | 1 | 1 |

That is, the investor can remove all funds from the three risky assets and convert it into cash (the risk-free asset), or the investor can take funds out of asset 1 and asset 2 and put it in asset 3, or the investor can take funds out from asset 1 and asset 3 and put it in asset 2, or the investor can take funds out from asset 2 and asset 3 and put it in asset 1, or the investor can take funds out from asset 1 and put it in asset 2 and asset 3, or the investor can take funds out of asset 2 and put it in asset 1 and asset 3, or the investor can take funds out of asset 3 and put it in asset 1 and asset 2, or the investor can convert some cash into funds and put it in all three of the assets. Thus there are 8 alternative choices for the investor when it comes to re-balancing his portfolio.

Of course, according to the Black-Scholes hedging principle, the re-balancing needs to be done each time by setting the optimal proportion of funds to be invested in each asset equal to the partial derivatives of the option valuation formula w.r.t. each of these assets. However, the total number of alternative choices available to the investor increases with every new risky asset that is added to the basket thereby contributing to the extrinsic utility in terms of the expected information content of the total portfolio.

**Coding of product information about multi-asset, structured financial portfolios**

Extending the entropy measure of extrinsic utility, we may conceptualize the interaction between the buyer and the vendor as a two-way communication flow whereby the vendor *informs* the buyer about the expected utility derivable from the product on offer and the buyer *informs* the seller about his or her individual expected utility criteria. An economic transaction goes through if the two sets of information are compatible. Of course, the greater expected information content of the vendor's communication, the higher is the extrinsic utility of the buyer. Intuitively, the expected information content of the vendor's communication will increase with increase in the variety of the product on offer, as that will increase the likelihood of matching the buyer's expected utility criteria.

The product information from vendor to potential buyer may be transferred through some medium e.g. the vendor's website on the Internet, a targeted e-mail or a telephonic promotion scheme. But such transmission of information is subject to noise and distractions brought about by environmental as well as psycho-cognitive factors. While a **distraction** is prima facie predictable, (e.g. the pop-up windows that keep on opening when some commercial websites are accessed), *noise* involves unpredictable perturbations (e.g. conflicting product information received from any competing sources).

Transmission of information calls for some kind of *coding*. Coding may be defined as a mapping of words from a source alphabet A to a code alphabet B. A **discrete, finite memory-less channel** with finite inputs and output alphabets is defined by a set of transition probabilities $p_i(j)$, i = 1, 2 … a and j = 1,2 … b with $\Sigma_j p_i(j) = 1$ and $p_i(j) \geq 0$. Here $p_i(j)$ is the probability that for an input letter i output letter j will be received.

A **code word** of length n is defined as a sequence of n input letters which are actually n integers chosen from 1,2 … a. A **block code** of length n having M words is a mapping of the message integers from 1 to M into a set of code words each having a fixed length n. Thus for a structured product with N component assets, a block code of length n having N words would be used to map message integers from 1 to N, corresponding to each of the N assets, into a set of a fixed-length code words. Then there would be a total number of $C = 2^N$ possible combinations such that $\log_2 C = N$ binary-state devises (flip-flops) would be needed.

A **decoding system** for a block code is the inverse mapping of all output words of length n into the original message integers from 1 to M. Assuming all message integers are used with same probability 1/M, the probability of error $P_e$ for a code and decoding system ensemble is defined as the probability of an integer being transmitted and received as a word which is mapped into another integer i.e. $P_e$ is the probability of wrongly decoding a message.

Therefore, in terms of our structured product set up, $P_e$ might be construed as the probability of misclassifying the best performing asset. Say within a structured product consisting of three risky assets - a blue-chip equity portfolio, a market-neutral hedge fund and a commodity future (and cash as the risk-free asset), while the original transmitted information indicates the hedge fund to be the best performer, due to erroneous decoding of the encoded message, the equity portfolio is interpreted as the best performer. Such erroneous decoding could result in investment funds being allocated to the wrong asset at the wrong time.

**The relevance of Shannon-Fano coding to product information transmission**

By the well-known Kraft's inequality we have **K = $\Sigma_n$ 2$^{-l_i}$ ≤ 1**, where $l_i$ stands for some definite code word lengths with a radix of 2 for binary encoding. For block codes, $l_i$ = l for i = 1, 2 … n. As per ***Shannon's coding theorem***, *it is possible to encode all sequences of n message integers into sequences of binary digits in such a way that the average number of binary digits per message symbol is approximately equally to the entropy of the source, the approximation increasing in accuracy with increase in n*. For efficient binary codes, K = 1 i.e. log$_2$ K = 0 as it corresponds to the maximal entropy condition. Therefore the inequality occurs if and only if **$p_i$ ≠ 2$^{-l_i}$**. Though the Shannon-Fano coding scheme is not strictly the most efficient, it has the advantage of directly deriving the code word length $l_i$ from the corresponding probability $p_i$. With source symbols $s_1$, $s_2$ … $s_n$ and their corresponding probabilities $p_1$, $p_2$ … $p_n$, where for each $p_i$ there is an integer $l_i$, then given that we have bounds that span an unit length, we have the following relationship:

$$\log_2 (p_i^{-1}) \leq l_i < \log_2 (p_i^{-1}) + 1 \tag{15}$$

Removing the logs, taking reciprocals and summing each term we therefore get,

$$\Sigma_n\, p_i \geq \Sigma_n\, 2^{l_i} \geq p_i/2, \text{ that is,}$$

$$1 \geq \Sigma_n\, 2^{l_i} \geq \tfrac{1}{2} \tag{16}$$

Inequality (16) gets us back to the Kraft's inequality. This shows that there is an instantaneously decodable code having the Shannon-Fano lengths $l_i$. By multiplying inequality (15) by $p_i$ and summing we get:

$$\Sigma_n\, (p_i \log_2 p_i^{-1}) \leq \Sigma_n\, p_i l_i < \Sigma_n\, (p_i \log_2 p_i^{-1}) + 1, \text{ that is,}$$

$$H_2(S) \leq L \leq H_2(S) + 1 \tag{17}$$

That is, in terms of the average Shannon-Fano code length L, we have conditional entropy as an effective lower bound while it is also the non-integral component of the upper bound of L. This underlines the relevance of a Shannon-Fano form of coding to our structured product formulation as this implies that the average code word length used in this form of product information coding would be **bounded by a measure of extrinsic utility** to the potential investor of the structured financial product itself, which is definitely an intuitively appealing prospect.

**Conceptualizing product information transmission as a Markov process**

The Black-Scholes option-pricing model is based on the underlying assumption that asset prices evolve according to the geometric diffusion process of a Brownian motion. The Brownian motion model has the following fundamental assumptions:

(1). $W_0 = 0$

(2). $W_t - W_s$ is a random variable that is normally distributed with mean 0 and variance t-s

(3). $W_t - W_s$ is independent of $W_v - W_u$ if (s, t) and (u, v) are non-overlapping time intervals

Property (3) implies that the Brownian motion is a *Markovian process* with no long-term memory. The switching behavior of asset prices from "high" (Bull state) to "low" (Bear state) and vice versa according to Markovian transition rule constitutes a well-researched topic in stochastic finance. It has in fact been proved that a steady-state equilibrium exists when the state probabilities are equalized for a stationary transition-probability matrix (Bhattacharya, 2001). This steady-state equilibrium corresponds to the condition of strong efficiency in the financial markets whereby no historical market information can result in arbitrage opportunities over any significant length of time.

By logical extension, considering a structured portfolio with n assets, the best performer may be hypothesized to be traceable by a first-order Markov process, whereby the best performing asset at time t+1 is dependent on the best performing asset at time t. For example, with n = 3, we have the following state-transition matrix:

|        | Asset 1   | Asset 2   | Asset 3   |
|--------|-----------|-----------|-----------|
| Asset 1 | P (1 \| 1) | P (2 \| 1) | P (3 \| 1) |
| Asset 2 | P (2 \| 1) | P (2 \| 2) | P (3 \| 2) |
| Asset 3 | P (3 \| 1) | P (3 \| 2) | P (3 \| 3) |

In information theory also, a similar Markov structure is used to improve the encoding of a source alphabet. For each state in the Markov system, an appropriate code can be obtained from the corresponding transition probabilities of leaving that state. The efficiency gain will depend on how variable the probabilities are for each state. However, as the order of the Markov process is increased, the gain will tend to be less and less while the number of attainable states approach infinity.

The strength of the Markov formulation lies in its capacity of handling correlation between successive states. If $S_1, S_2 \ldots S_m$ are the first m states of a stochastic variable, what is the probability that the next state will be $S_i$? This is written as the conditional probability $p(S_i | S_1, S_2 \ldots S_m)$. Then, the Shannon measure of information from a state $S_i$ is given as usual as follows:

$$I(S_i | S_1, S_2 \ldots S_m) = \log_2 \{p(S_i | S_1, S_2 \ldots S_m)\}^{-1} \qquad (17)$$

The entropy of a Markov process is then derived as follows:

$$H(S) = \sum_{S^{m+1}} p(S_1, S_2 \ldots S_m, S_i) I(S_i | S_1, S_2 \ldots S_m) \qquad (18)$$

Then the extrinsic utility to an investor from a structured financial product expressed in terms of the entropy of a Markov process governing the state-transition of the best performing asset over N component risky assets (and cash as the one risk-free asset) within the structured portfolio would be given as follows:

$$U_x = H(\text{Portfolio}) = \sum_{S^{N+1}} p(S_1, S_2 \ldots S_m, S_i) I(S_i | S_1, S_2 \ldots S_m) \qquad (19)$$

However, to find the entropy of a Markov source alphabet one needs to explicitly derive the stationary probabilities of being in each state of the Markov process. But these state probabilities may be hard to derive explicitly especially if there are a large number of allowable states (e.g. corresponding to a large number of elementary risky assets within a structured financial product). Using **Gibbs inequality**, it can be show that the following limit can be imposed for bounding the entropy of the Markov process:

$\Sigma_j$ **p ($S_j$) H (Portfolio | $S_j$) $\leq$ H ($S^*$)**, where **H ($S^*$)** is termed the ***adjoint system***   (20)

The entropy of the original message symbols given by the zero memory source adjoint system with **p ($S_i$) = $p_i$** bound the entropy of the Markov process. The equality holds if and only if **p ($S_j$, $S_i$) = $p_j p_i$** that is, in terms of the structured portfolio set up, the equality holds if and only if the joint probability of the best performer being the pair of assets i and j is equal to the product of their individual probabilities (Hamming, 1986). Thus a clear analogical parallel may be drawn between Markovian structure of the coding process and performances of financial assets contained within a structured investment portfolio.

## Conclusion and scope for future research

In this paper we have basically outlined a novel methodological approach whereby expected information measure is used as a measure of utility derivable from a basket of commodities. We have illustrated the concepts with an applied finance perspective whereby we have used this methodological approach to derive a measure of investor utility from a structured financial portfolio consisting of many elementary risky assets combined with cash as the risk-free asset thereby giving the product a quasi - capital guarantee status. We have also borrowed concepts from mathematical information theory and coding to draw analogical parallels with the utility structures evolving out of multi-asset, structured financial products. In particular, principles of Shannon-Fano coding have been applied to the coding of product information for transmission from vendor (fund manager) to the potential buyer (investor). Finally we have dwelled upon the very similar Markovian structure of coding process and that of asset performances.

This paper in many ways is a curtain raiser on the different ways in which tools and concepts from mathematical information theory can be applied in utility analysis in general and to analyzing investor utility preferences in particular. It seeks to extend the normal peripheries of utility theory to a new domain – that of information theoretic utility. Thus a cardinal measure of utility is proposed in the form of the Shannon-Boltzmann entropy measure. Being a new methodological approach, the scope of future research is boundless especially in exploring the analogical Markovian properties of asset performances and message transmission and devising an efficient coding scheme to represent the two-way transfer of utility information from vendor to buyer and vice versa. The mathematical kinship between neoclassical utility theory and classical thermodynamics is also worth exploring, may be aimed at establishing some higher-dimensional, theoretical connectivity between the isotherms and the indifference curves!